\documentstyle[11pt]{article}
\def\hang{\hangindent\parindent}
\def\textindent#1{\indent\llap{#1\enspace}\ignorespaces}

\begin{document}

\centerline{\huge\bf Global existence and exponential decay
  } \vskip 5pt \centerline{\huge\bf of the solution for
a viscoelastic wave } \vskip 5pt \centerline{\huge\bf equation with
a delay
}

\vspace{0.5cm} \centerline{  Qiuyi Dai$^\dag$\ \ \ \ Zhifeng
Yang$^{\dag\ \ddag}$ \footnote[1]{\scriptsize Corresponding author.
E-mail: zhifeng \underline \ \ yang@126.com.}
}

\vskip 3pt \centerline{\footnotesize\it $\dag$ College of
Mathematics and Computer Science, Hunan Normal University}
\centerline{\footnotesize\it Changsha, Hunan, 410081, P.R.China}
\centerline{\footnotesize\it $\ddag$ Department of Mathematics and
Computational Science, Hengyang Normal University}
\centerline{\footnotesize\it Hengyang, Hunan, 421008, P.R.China}

\vspace{0.5cm}

\begin{quote}
{\bf Abstract}~~ {\small  In this paper, we consider initial-boundary value problem of viscoelastic wave equation
with a delay term in the interior feedback. Namely, we study the
following equation
$$u_{tt}(x,t)-\Delta u(x,t)+\int_0^t g(t-s)\Delta u(x,s)ds +\mu_1 u_t(x,t)+ \mu_2 u_t(x,t-\tau)=0$$
together with initial-boundary conditions of Dirichlet type in
$\Omega\times (0,+\infty)$, and prove that for arbitrary real
numbers $\mu_1$ and $\mu_2$, the above mentioned problem has a
unique global solution under suitable assumptions on the kernel $g$.
This improve the results of the previous literature such as [6] and
[13] by removing the restriction imposed on $\mu_1$ and $\mu_2$.
Furthermore, we also get an exponential decay results for the energy
of the concerned problem in the case $\mu_1=0$ which solves an open
problem proposed by M. Kirane and B. Said-Houari in [13]. }



\end{quote}

\vspace{0.3cm}

\section{Introduction}

~~~~It is well-known that the free vibration of membrane can be
described by wave equation of the form
$$u_{tt}(x,t)-\Delta u(x,t)=0 \eqno(1.1)$$
in $\Omega \times (0,+\infty)$, subjected to initial conditions and
some boundary conditions. Here, $\Omega \subset {\mathbf{R}}^2$ is a
bounded domain. For the case of high dimension, there are also many
discussions in mathematical physics. For the better comprehension of
our motivation, we appeal to readers to keep in mind that the system
(1.1) is conservative.

To study the propagation mechanism of wave, two main factors, that
is damping and time-delay, are taken into consideration. When the
damping occurs, the corresponding mathematical model becomes
$$u_{tt}(x,t)-\Delta u(x,t)+\mu u_t=0, \eqno(1.2)$$
and the energy of this system is dissipative. We here point out that
the results about the models (1.1) and (1.2) are quite abundant. And
what is regretful is that we cannot list these literatures one by
one.

Moreover, time-delay effects often appears in many daily life
practical problems. These hereditary effects are sometime
unavoidable, and it may turn a well-behaved system
into a wild one. That is to say, they can induce some
instabilities(see for instance[1-5] and the references therein). The
stability issue of systems with delay is, therefore, of theoretical
and practical importance.

For the wave equation with boundary or interior delay, S. Nicaise,
C. Pignotti, M. Gugat and G. Leugering et al. obtained some profound
results[6-8,20-22]. In [6], the authors examined a system of wave
equation with the damping and the time-delay inside the domain.
Namely, they considered the following system
$$
\left\{
\begin{array}{l}
u_{tt}(x,t)-\Delta u(x,t)+\mu_1 u_t(x,t)+ \mu_2 u_t(x,t-\tau)=0,\ x\in \Omega,\ t>0,\nonumber \\
u(x,t)=0,\ \ \ \ \ \ \ \ \ \ \ \ \ \ \ \ \ \ \ \ \ \ \ \ \ \ \ \ \ \ \ \ \ \ \ \ \ \ \ \ \ \ \ \ \ \  \ \ \ \ \ \ \ \ x\in \partial \Omega,\ t\geq 0,\nonumber \\
u(x,0)=u_{0}(x),\ \ u_t(x,0)=u_{1}(x),\ \ \ \ \ \ \ \ \ \ \ \ \ \ \
\ \ \
 \ \ \ \ \ \ x\in \Omega ,\nonumber \\
u_{t}(x,t-\tau)=f_0(x,t-\tau),\ \ \ \ \ \ \ \ \ \ \ \ \ \ \ \ \ \ \
\ \ \ \ \ \ \ \ \ \ \ \ \ \ \ \ x\in \Omega ,\ 0<t<\tau, \nonumber
\end{array}
\right.\eqno(1.3)
$$
and proved that the energy of the above problem is exponentially
decay (or the trivial solution is exponentially stable) provided
that $0<\mu_2<\mu_1$. Furthermore, it is also showed in the case
$\mu_2 \geq\mu_1$ that there exists a sequence of arbitrary small
(or large) delays such that instabilities occur. The same results
were showed when both the damping and the delay act on the boundary.
Particularly, if $\mu_1 =0$, that is, if we have only the delay part
in the boundary condition, the system becomes unstable (see [2]).
From these results, we may figure out that the delay term lead to
instability of the system. And to avoid this problem, M. Gugat
considered feedback laws where a certain delay is included as a part
of the control law and not as a perturbation, and obtained the
exponential stability of the proposed feedback with retarded input
(see [20]). In addition, M. Gugat and G. Leugering considered the
feedback stabilization of quasilinear hyperbolic systems on
star-shaped networks. They obtained the exponential stability of the
system by introducing an $L^2$-Lyapunov function with delay terms
(see [21]). Recently, C. Pignotti considered the wave equation with
internal distributed time delay and local damping in a bounded and
smooth domain, and showed that an exponential stability result holds
if the coefficient of the delay term is sufficiently small(see
[22]).

In this paper, we are concerned with the initial-boundary value
problem(IBVP) of the viscoelastic wave equation
$$
\left\{
\begin{array}{l}
u_{tt}(x,t)-\Delta u(x,t)+\int_0^t g(t-s)\Delta u(x,s)ds\nonumber \\
\ \ \ \ \ \ \ \ \ \ \ +\mu_1 u_t(x,t)+ \mu_2 u_t(x,t-\tau)=0,\ \ \ \ \ x\in \Omega,\ t>0,\nonumber \\
u(x,t)=0,\ \ \ \ \ \ \ \ \ \ \ \ \ \ \ \ \ \ \ \  \ \ \ \ \ \ \ \ \ \  \ \ \ \ \ \ \ \ \ \ \ \ \ \  x\in \partial \Omega,\ t\geq 0,\nonumber \\
u(x,0)=u_{0}(x),\ \ u_t(x,0)=u_{1}(x),\ \ \ \  \ \ \ \ \ \ \ \ \ \   x\in \Omega ,\nonumber \\
u_{t}(x,t-\tau)=f_0(x,t-\tau),\ \ \  \ \ \ \ \ \ \ \ \ \  \ \ \ \ \
\ \ \ \ \  \ x\in \Omega ,\ 0<t<\tau, \nonumber
\end{array}
\right.\eqno(1.5)
$$
where $\Omega\subset \textbf{R}^n, n\geq 1$ be a regular and bounded
domain with a  boundary $\partial\Omega $ of class $C^2$. Moreover,
$\mu_1,\mu_2$ are two real numbers, $\tau>0$ represents the time
delay and $u_0,u_1,f_0$ are given initial data belonging to suitable
spaces.

System (1.5) comes from linear models for propagation of
viscoelastic wave in a compressible fluid. In mechanics, it is well
known that solid and fluid materials exhibit not only elasticity but
also hereditary properties. The hereditary properties are described
by viscoelasticity, where the mechanical response of the materials
is taken to be influenced by the previous behavior of the materials
themselves. In other words, the viscoelastic materials exhibit
memory effects. This leads to a constitutive relationship between
the stress and strain involving the convolution of the strain with a
relaxation function. This term is called a ``memory term". From the
mathematical point of view, these hereditary properties are modeled
by integro-differential operators.

To motivate our present work, let us recall some previous results regarding the viscoelastic wave equation. In the absence of time-delay,
the viscoelastic wave equation of the form
$$u_{tt}-\Delta u+\int_0^t g(t-s)\Delta u(x,s)ds+a(x)u_t=0, \eqno(1.6)$$
in $\Omega\times(0,\infty)$, subjected to initial conditions and
boundary conditions of Dirichlet type has been considered by
Cavalcanti et al.[9]. They showed an exponential decay result under
some restrictions on $a(x)$ and $g(t)$. To be specific, they assumed
that $a: \Omega\rightarrow \mathbf{R}$ is a nonnegative and bounded
function, and the kernel $g$ in the memory term decays
exponentially. This result has been improved later by Berrimi and
Messaoudi[10] under weaker conditions on both $a$ and $g$. And then,
more general problems than the one considered in [9] are studied in
[11,12]. Some profound results about the relationship between the
stability of solutions and the relaxation kernel $g$ are obtained.
It is worthwhile to note that, if $a(x)\equiv 0$ in (1.6), some
results about the energy decay are obtained by many researchers.
See, for instance Tatar[18], Mustafa and Messaoudi[19].

Then a nature problem is that what would happen when a delay term
occurs in (1.6). In this case, M. Kirane and B. Said-Houari studied
problem (1.5) in [13] with coefficients $\mu_1$ and $\mu_2$
positive. They showed that problem (1.5) has a unique weak solution
(see [13],Theorem 3.1), and its energy is exponentially decay
provided that $\mu_2\leq \mu_1$.
 It is worth pointing out that the assumptions $\mu_1>0$ and $\mu_2\leq \mu_1$ play a decisive role in proof of the above mentioned results.
 Because the authors must use the damping term $\mu_1 u_t(x,t) $ to control the delay
 term $\mu_2u_t(x,t-\tau)$ in the priori estimate of the solution and the decay estimate of the energy.
 This makes the authors arise an open problem that whether the decay property of the energy they
have obtained are preserved in the case $\mu_1=0$? By the way, the
results of [13] has been generalized recently by W.J. Liu to a
system with time dependent delay in [14].
 The method used by W.J. Liu is very similar to those used in [13].

The aims of the present paper are twofold. At first, we aim to prove
an existence result of problem (1.5) without restrictions of
$\mu_1,\mu_2>0$ and $\mu_2\leq \mu_1$. Unlike in [13], our new idea
is to control the delay term by the derivative of the energy instead
of by the damping term in the priori estimate of solutions. In the
second, we will give a positive answer to the open problem proposed
by M. Kirane and B. Said-Houari in [13]. That is we will prove a
energy decay result for problem (1.5) in the case $\mu_1=0$. The
main difficulty in handling this problem is that we have no damping
term to control the delay term in the estimate of the energy decay.
To overcome this difficulty, our basic idea is to control the delay
term by making use of the viscoelasticity term. And in order to
achieve this goal, a restriction of the size between the parameter
$\mu_2$ and the kernel $g$ and a new Lyapunov functional, different
from that one used in [13], are needed.

The paper is organized as follows. In Section 2 we give some
preliminaries. In Section 3, we prove the existence and the
uniqueness of the solution of problem (1,5). In Section 4, we prove
the exponential decay of the energy of this problem in the case
$\mu_1=0$.

\bigbreak
\section{Preliminaries}
\bigbreak

~~~~In this section, we shall prepare some materials needed in the
proof of our result and state. We use the standard Lebesgue space
$L^2(\Omega)$ and Sobolev space $H^1_0(\Omega)$ with their usual
norms $\|\cdot\|_2$ and $\|\cdot\|_{H_0^1}$ respectively. And we
will write $(\cdot,\cdot)$ to denote the inner product in
$L^2(\Omega)$.
Throughout this paper, $C$ and $C_i$ are used to denote generic
positive constants.

We first state the general assumptions on the kernel function $g$ as
follows:

(A1) $g:\textbf{R}^+\rightarrow \textbf{R}^+$ is a $C^1$ function
satisfying
$$g(0)>0,\ \ \ \ 1-\int_0^\infty g(s)ds=l>0.$$

(A2) There exists a positive constant $\zeta$ such that
$$g^\prime(t)\leq -\zeta g(t), \ \ \forall t>0. \eqno(2.1)$$

A typical example of such function is
$$g(t)=e^{-at}, a>1.$$

Now,let us introduce the following notation:
$$(\phi\circ\psi)(t):=\int_0^t \phi(t-s)\int_\Omega |\psi(t)-\psi(s)|^2dx ds.$$

For any function $\phi\in C^1(\mathbf{R})$ and any $\psi \in
H^1(0,T,L^2(\Omega))$, straightforward computations lead to(see
[13,15])
$$\int_\Omega\psi_t(t)\int_0^t \phi(t-s) \psi(s)dsdx\ \ \ \ \ \ \ \ \ \ \ \ \ \ \ \ \ \ \ \ \ \ \ \ \ \ \ \ \ \ \ \ \ \ \ \ \ \ \ \ \ \ \ \ \ \ \ \ \ \ \ \ \ \ \ \ \ \ \ \ \ $$
$$=\frac{1}{2}(\phi^\prime\circ\psi)(t)-\frac{1}{2}\phi(t)\|\psi\|_2^2-\frac{1}{2}\frac{d}{dt}\left\{(\phi\circ\psi)(t)-\left(\int_0^t \phi(s)ds\right)\|\psi\|_2^2\right\}.\eqno(2.2)$$

The following lemma will be used later in order to define the new
modified energy functional of problem (1.5).



{\bf Lemma 2.1}(see [13, Lemma 2.2] or [16]). For $u\in
L^2(0,T,H_0^1(\Omega))$, we have
$$\int_\Omega\left(\int_0^t g(t-s)(u(t)-u(s))ds\right)^2dx\leq (1-l)C_\ast^2(g\circ \nabla u)(t),\eqno(2.3)$$
where $C_\ast$ is the Poincar\'{e} constant and $l$ is given in
(A1).

Let
$$U= \left\{u|u\in L^2\left(0,T;H_0^1(\Omega)\right), \ \ u_t\in
L^2\left(0,T;L^2(\Omega)\right)\right\}$$ and
$$V=\left\{v|v\in U,\ \  v(x,T)=0,\ \  v_t(x,T)=0 \right\}.$$

Next, we give the definition of the weak solution of problem (1.5).

{\bf Definition 2.2} We say a function $u\in U $
is a weak solution of the
initial-boundary value problem (1.5) provided

$$-\int_0^T(  u_{t}, v_t)dt+\int_0^T(\nabla u(t),\nabla v(t))dt-\int_0^T\int_0^t
g(t-s)(\nabla u(s),\nabla v(s))ds dt$$
$$\ \ \ \ \ \ \ \ \ \ \ \ \ \ \ \ \ \ \ \ \ \ \ \ \ \ \ \ +\mu_1\int_0^T(u_t(t),v)dt+\mu_2\int_0^T(u_t(t-\tau),v)dt+(v(x,0),u_1(x))=0 \eqno(2.4)$$
for each $v\in V$ and
$$u(x,0)=u_0,\ \ u_t(x,0)=u_1,\ \ u_t(x,t-\tau)=f_0(x,t-\tau),\ \  t\in
(0,\tau).\eqno(2.5)$$

{\emph{\textbf{Remark 2.3}}}~~\emph{In view of Theorem 2 in $\S5.92$
in [23], we know $u\in C\left(0,T;L^2(\Omega)\right)$ and $u_t\in
C\left(0,T;H^{-1}(\Omega)\right)$. Consequently the equalities (2.5)
above make sense. }

\bigbreak
\section{Existence of weak solutions }
\bigbreak

~~~In this section, we will give a sufficient condition that
guarantees the well-posedness of the problem (1.5). We will use the
Faedo-Galerkin approximation to prove it. So we will need the
following prior estimates.

{\bf Theorem 3.1}(Prior estimates) For all $u\in U$, there exists a
constant $C$, depending only on $\Omega$ and $T$, such that
$${\mathcal{E}}(u(t))\leq C,$$
where
$${\mathcal{E}}(u(t))=\frac{1}{2}\left(\|u_{t}\|_2^2
+\left(1-\int_0^t g(s)ds\right)\|\nabla u\|_2^2 +(g\circ \nabla
u)(t)\right).\eqno(3.1)$$

{\bf \emph{\textbf{Proof}}}~~Since functions of the space
$C^\infty((0,T)\times\Omega)$ are dense in the space $U$, we just
need to prove the conclusion in the case of $u\in
C^\infty((0,T)\times\Omega)$. So, the following calculations are all
in this space.

Multiplying the first equation in (1.5) by $u_t$ and integrating on
$\Omega$, straightforward computations lead to
$${\mathcal{E}}^\prime (u(t))=\frac{1}{2}\left((g^\prime\circ \nabla u)(t)-g(t)\|\nabla u(t)\|_2^2 \right)-\mu_1\int_\Omega u^2_{t}(t)dx-\mu_2\int_\Omega u_{t}(t-\tau)u_{t}(t)dx.\eqno(3.2)$$
By the assumptions (A1) and (A2) about $g(t)$, we have
$$ (g^\prime\circ
\nabla u)(t)- g(t)\|\nabla u(t)\|_2^2 <0.\eqno(3.3)$$
 So, thanks to the Cauchy-Schwartz inequality, we deduce
$${\mathcal{E}}^\prime (u(t))\leq \frac{|\mu_2|}{2}\int_\Omega
u^2_{t}(t-\tau)dx+\left(\frac{|\mu_2|}{2}+|\mu_1|\right)\int_\Omega
u^2_{t}(t)dx.\eqno(3.4)$$

By integrating on $[0,t]$, we have
$${\mathcal{E}}(u(t))\leq \frac{|\mu_2|}{2}\int_0^t \int_\Omega
u^2_{s}(s-\tau)dxds+\left(\frac{|\mu_2|}{2}+|\mu_1|\right)\int_0^t\int_\Omega
u^2_{s}(s)dxds+{\mathcal{E}}(u(0)),\eqno(3.5)$$ where
${\mathcal{E}}(u(0))$ is a nonnegative constant.

Now, using the past history values about $u_{t}(t),\ \ t\in
[-\tau,0]$, the first term in the right-hand side of (3.5) can be
rewritten as follows
$$\int_0^t \int_\Omega u^2_{s}(s-\tau)dxds\ \ \ \ \ \ \ \ \ \ \ \ \ \ \ \ \ \ \ \ \ \ \ \ \ \ \ \ \ \ $$
$$=\int_\Omega \int_{-\tau}^{t-\tau}u^2_{\rho}(\rho)d\rho dx\ \ \ \ \ \ \ \ \ \ \ \ \ \ \ \ \ \ \ \ \ \ \ \ \ \ \ $$
$$=\int_\Omega \int_{-\tau}^{0}u^2_{\rho}(\rho)d\rho dx+\int_\Omega \int_{0}^{t-\tau}u^2_{\rho}(\rho)d\rho dx$$
$$=\int_\Omega \int_{-\tau}^{0}f^2_{0}(\rho)d\rho dx+\int_\Omega \int_{0}^{t-\tau}u^2_{\rho}(\rho)d\rho dx$$
$$\leq \int_\Omega \int_{-\tau}^{0}f^2_{0}(\rho)d\rho dx+\int_\Omega \int_{0}^{t}u^2_{\rho}(\rho)d\rho dx.\ \ \   \eqno(3.6)$$

 From (3.5) and (3.6), we obtain
$${\mathcal{E}}(u(t))\leq (|\mu_2|+|\mu_1|)\int_\Omega \int_{0}^{t}u^2_{\rho}(\rho)d\rho dx+\frac{|\mu_2|}{2}\int_\Omega \int_{-\tau}^{0}f^2_{0}(\rho)d\rho dx+{\mathcal{E}}(u(0)).\eqno(3.7)$$
Thus, we have
$${\mathcal{E}}(u(t))\leq
2(|\mu_2|+|\mu_1|)\int_0^t{\mathcal{E}}(u(s))ds+\frac{|\mu_2|}{2}\int_\Omega
\int_{-\tau}^{0}f^2_{0}(\rho)d\rho
dx+{\mathcal{E}}(u(0)).\eqno(3.8)$$
 By the Gronwall
inequality, once $T>0$ be given, $\forall t \in [0,T]$, we have
$${\mathcal{E}}(u(t))\leq \left(\frac{|\mu_2|}{2}\int_\Omega
\int_{-\tau}^{0}f^2_{0}(\rho)d\rho
dx+{\mathcal{E}}(u(0))\right)e^{2(|\mu_2|+|\mu_1|)T},\ \ \forall
t\in [0,T].\eqno(3.9)$$ Denote
$$C=\left(\frac{|\mu_2|}{2}\int_\Omega
\int_{-\tau}^{0}f^2_{0}(\rho)d\rho
dx+{\mathcal{E}}(u(0))\right)e^{2(|\mu_2|+|\mu_1|)T}.$$ Then we have
${\mathcal{E}}(u(t))\leq C.$ $\Box$

With the above prior estimates, it is easy to prove the uniqueness
of solutions as follows:

{\bf Theorem 3.2}(uniqueness of solution). There exists at most one
solution of problem (1.5) in the sense of \emph{definition 2.2}.

{\bf \emph{\textbf{Proof}}}~~It suffices to show that the only weak
solution of (1.5) with $u_0=u_1=f_0=0$ is
$$u\equiv 0.\eqno(3.10)$$


According to the energy estimate (3.9) in Theorem 3.1, and noting
that $f_0=0,\ \  {\mathcal{E}}(u(0))=0$, we obtain
$${\mathcal{E}}(u(t))=0, \ \ \forall t\in [0,T].\eqno(3.11)$$
So, we have
$$\|u_t(t)\|_2=\|\nabla u(t)\|_2=0, \ \ \forall t\in [0,T].\eqno(3.12)$$
And this implies (3.10). Thus, we conclude that problem (1.5) has at
most one solution. $\Box$

Next, we give the existence of weak solution by the following
theorem.

 {\bf Theorem 3.3}(Existence of weak solution). Assume that (A1) and (A2) hold. Then given $u_0\in H_0^1(\Omega)$, $u_1\in
L^2(\Omega)$, $f_0 \in L^2(-\tau,0;\Omega)$ and $T>0$, there exists
a weak solution $u$, defined by \emph{definition 2.2}, of problem
(1.5) on $(0,T)$.

{\bf \emph{\textbf{Proof}}}~~We first construct the approximate
solution of (1.5) by Galerkin's method. Let
$\{\lambda_j\}_{j=1}^\infty$ be the eigenvalue of the following
eigenvalue problem:
$$
\left\{
\begin{array}{l}
-\Delta w=\lambda w, \ \ \ x\in \Omega,\nonumber \\
w=0,\ \ \ \ \ \ \ \ \ \ x\in \partial\Omega.\nonumber
\end{array}
\right. \eqno(3.13)
$$
 and $w_j$ is the eigenfunctions associated with
 $\lambda_j$($j=1,2,\cdots$). Furthermore, we choose
 $\{w_j\}_{j=1}^\infty$ such that
$$
\int_\Omega w_jw_k dx= \left\{
\begin{array}{l}
1, \ \ \ j=k,\nonumber \\
0,\ \ \  j\neq k.\nonumber
\end{array}
\right.\eqno(3.14)
$$
So, we have
$$
\int_\Omega \nabla w_j\nabla w_k dx= \left\{
\begin{array}{l}
\lambda_j, \ \ \ j=k,\nonumber \\
0,\ \ \ \  j\neq k.\nonumber
\end{array}
\right.\eqno(3.15)
$$
Moreover, $\{w_j\}_{j=1}^\infty$ is the orthogonal bases in
$L^2(\Omega)$ and in $H_0^1(\Omega)$.

We will seek an approximate solution in the form
$$u_n(t)=\sum_{j=1}^n \gamma_n^j (t)w_j, \eqno(3.16)$$
where we intend to select the coefficients $\gamma_n^j(t)$ ($0\leq
t\leq T, j=1,\cdots,n$) to satisfy
$$\gamma_n^j(0)=(u_0,w_j),\ \ \ \  \gamma_{nt}^j(0)=(u_1,w_j),\eqno(3.17)$$
$$ \gamma_{nt}^j(x,t-\tau)=f_0(x,t-\tau),\ \ \ \ t\in (0,\tau), \eqno(3.18)$$
and
$$( u_{ntt}, w_j)+(\nabla u_n(t),\nabla w_j(t))-\int_0^t
g(t-s)(\nabla u_n(s),\nabla
w_j)ds+\mu_1(u_{nt}(t),w_j)+\mu_2(u_{nt}(t-\tau),w_j)=0.
\eqno(3.19)$$

By (3.16), (3.19) becomes the linear system of ODE
$$\ddot{\gamma}_{n}^{j}(t)+\mu_1 \dot{\gamma}_n^j(t)+\mu_2 \dot{\gamma}_n^j(t-\tau)+\lambda_j\left(1-\int_0^t g(s)ds\right)\gamma_n^j=0,\eqno(3.20)$$
subject to the initial conditions (3.17) and the history value
(3.18).

For the initial value problem (3.20), (3.17) and (3.18), we can use
the contraction mapping principle to prove that there exists a
unique solution
$\gamma_n(t)=\left(\gamma_n^1(t),\cdots,\gamma_n^n(t)\right)\in C^2$
 when $t$ is small enough. And then, by the above prior estimates(see
Theorem 3.1), we have $ \| \gamma_n \|_{H_0^1}\leq C$. So, we can
extend this solution to the whole interval $[0,T]$.


Next, according to the prior estimates, we see that the sequence
$\{u_n\}_{n=1}^\infty$ is bounded in
$L^2\left(0,T;H_0^1(\Omega)\right)$, $\{u_{nt}\}_{n=1}^\infty$ is
bounded in $L^2\left(0,T;L^2(\Omega)\right)$.

As a consequence there exists a subsequence $\{u_k\}_{k=1}^\infty
\subset \{u_n\}_{n=1}^\infty$ and $u\in
L^2\left(0,T;H_0^1(\Omega)\right)$, with $u_t\in
L^2\left(0,T;L^2(\Omega)\right)$ such that
$$u_k\rightharpoonup u\ \ weakly\ \ in\ \ L^2\left(0,T;H_0^1(\Omega)\right),\ \ as\ \ k\rightarrow \infty,\eqno(3.21)$$
$$u_{kt}\rightharpoonup u_t\ \ weakly\ \ in\ \ L^2\left(0,T;L^2(\Omega)\right),\ \ as\ \ k\rightarrow \infty.\eqno(3.22)$$

Next, arguing as in [17, Theorem 3.1], we will prove that $u$ is the
weak solution of problem (1.5). For this, we choose a function
$v\in\{\eta|\eta\in
C^1\left(0,T;H_0^1(\Omega)\right),\eta(x,T)=0,\eta_t(x,T)=0\}$ of
the form
$$v(t)=\sum_{j=1}^N \gamma^j(t)w_j,\eqno(3.23)$$
where $\{\gamma^j\}_{j=1}^N$ are smooth functions, and $N$ is a
fixed integer.

We select $k\geq N$, multiply the equation (3.19) by $\gamma^j(t)$,
sum $j=1,\cdots,N$, and then integrate with respect to $t$, to
discover
$$I+II=0, \eqno(3.24)$$
where
$$I=\int_0^T \left[-( u_{kt},
v_t)+(\nabla u_k(t),\nabla v(t))-\int_0^t g(t-s)(\nabla
u_k(s),\nabla v(s))ds\right]dt$$ and
$$II=\int_0^T \left[  \mu_1(u_{kt}(t),v)+\mu_2(u_{kt}(t-\tau),v) \right]dt+(v(x,0),u_{kt}(x,0)).$$


By (3.21),(3.22) and passing to the limit in (3.24) as $k\rightarrow
\infty$, we obtain
$$\int_0^T \left[-( u_{t},
v_t)+(\nabla u(t),\nabla v(t))-\int_0^t g(t-s)(\nabla u(s),\nabla
v(s))ds\right]dt$$
$$\ \ \ \ \ \ \ \ \ \ \ \ \ \ \ \ \ \ \ \ \ \ \ \ \ \ \ \ \ +\int_0^T \left[  \mu_1(u_{t}(t),v)+\mu_2(u_{t}(t-\tau),v) \right]dt+(v(x,0),u_{1}(x))=0 \eqno(3.25)$$
for all functions $v\in V$, since functions of the form (3.23) are
dense in this space. And thus, (2.4) holds.


We must now verify that the limit function $u$ satisfies the initial
conditions and the history value, i.e.
$$u(0)=u_0,\ \  u_t(0)=u_1,\eqno(3.26)$$
$$u_t(x,t-\tau)=f_0(x,t-\tau), \ \ t\in (0,\tau).\eqno(3.27)$$
For this, choose any function $v\in
C^2\left(0,T;H_0^1(\Omega)\right)$, with $v(T)=v_t(T)=0$. Then
integrating by parts with respect to $t$ in (3.25), we find
$$\int_0^T \left[(v_{tt},u)+(\nabla u(t),\nabla v(t))-\int_0^t g(t-s)(\nabla
u(s),\nabla v(s))ds\right]dt\ \ \ \ \ \ \ \ \ \ \ \ \ \ \ \ \ \ $$
$$+\int_0^T \left[  \mu_1(u_{t}(t),v)+\mu_2(u_{t}(t-\tau),v) \right]dt=(u(0),v_t(0))-( u_t(0),v(0)).\eqno(3.28)$$
Similarly from (3.24) we deduce
$$\int_0^T \left[(v_{tt},u_k)+(\nabla u_k(t),\nabla v(t))-\int_0^t g(t-s)(\nabla
u_k(s),\nabla v(s))ds\right]dt\ \ \ \ \ \ \ \ \ \ \ \ \ \ \ \ \ \ \
\ \ \ \ \ \ \ \ \ \ \ $$
$$+\int_0^T \left[  \mu_1(u_{kt}(t),v)+\mu_2(u_{kt}(t-\tau),v) \right]dt=(u_k(0),v_t(0))-(u_{kt}(0),v(0)).\eqno(3.29)$$
Recalling (3.17),(3.21) and (3.22), we get
$$\int_0^T \left[(v_{tt},u)+(\nabla u(t),\nabla v(t))-\int_0^t g(t-s)(\nabla
u(s),\nabla v(s))ds\right]dt\ \ \ \ \ \ \ \ \ \ \ \ \ \ \ \ \ \ $$
$$+\int_0^T \left[  \mu_1(u_{t}(t),v)+\mu_2(u_{t}(t-\tau),v) \right]dt=(u_0,v_t(0))-( u_1,v(0)),\eqno(3.30)$$
as $k\rightarrow \infty$. Comparing identities (3.28) and (3.30), we
conclude (3.26), since $v(0), v_t(0)$ are arbitrary. In addition,
$\forall j$, $\int_\Omega u_{kt}(t-\tau)w_j dx\rightharpoonup
\int_\Omega f_0(t-\tau)w_j dx$ weakly in $L^2(0,T;L^2(\Omega))$ as
$k\rightarrow \infty$. So we easily deduce
$u_t(x,t-\tau)=f_0(x,t-\tau)$. And then (3.27) holds. Hence $u$ is a
weak solution of (1.5) in the sense of \emph{Definition 2.2}. $\Box$


\bigbreak
\section{Exponential energy decay in the case $\mu_1=0$}
\bigbreak

~~~In this section, we will give the exponential energy decay result
about the problem (1.5) in the case $\mu_1=0$. That is to say, we
will study the following system
$$
\left\{
\begin{array}{l}
u_{tt}(x,t)-\Delta u(x,t)+\int_0^t g(t-s)\Delta u(x,s)ds + \mu u_t(x,t-\tau)=0,\ \  x\in \Omega,\ t>0,\nonumber \\
u(x,t)=0,\ \   x\in \partial \Omega,\ t\geq 0,\nonumber \\
u(x,0)=u_{0}(x),\ \ u_t(x,0)=u_{1}(x),\ \  x\in \Omega ,\nonumber \\
u_{t}(x,t-\tau)=f_0(x,t-\tau),\ \  x\in \Omega ,\ 0<t<\tau,
\nonumber
\end{array}
\right.\eqno(4.1)
$$
where $\mu$ is a real number.

We first define the classical energy by
$$e(t)=\frac{1}{2}\left(\|u_{t}\|_2^2+\|\nabla u\|_2^2\right),\ \  t\geq 0.$$
Then by the first equation in (4.1) it is easy to see that
$$e^\prime(t)=\int_\Omega \nabla u_t \int_0^t g(t-s) \nabla u(s) ds dx-\mu \int_\Omega u_t(t)u_t(t-\tau)dx,\ \ t \geq 0.$$
By (2.2) and inspired by [13], we define the new modified energy
functional of problem (4.1) as follows:
$$E(t)=\frac{1}{2}\left(\|u_{t}\|_2^2
+\left(1-\int_0^t g(s)ds\right)\|\nabla u\|_2^2 +(g\circ \nabla
u)(t)\right)+\frac{\xi}{2}\int_{t-\tau}^t\int_\Omega
e^{\sigma(s-t)}u_t^2(x,s) dx ds,\eqno(4.2)$$ where $\sigma$ and
$\xi$ are two positive constants to be determined later.

Since the function $g$ is positive, continuous and $g(0)>0$, then
for any $t\geq t_0>0$, we have
$$\int_0^t g(s)ds\geq \int_0^{t_0} g(s)ds=g_0.\eqno(4.3)$$

Our stability result reads as follows:

 {\bf Theorem 4.1}(Exponential energy decay) ~~Let $u$ be the solution of (4.1).
 Assume that $g$ satisfies (A1) and (A2), and $|\mu|< a$, where $a$ is a positive constant
defined by (4.31). Then, for any $t_0>0$, there exists two
 positive constants, $K$ and $\lambda$, such that the energy of
 problem (4.1) satisfies
 $$E(t)\leq Ke^{-\lambda \zeta (t-t_0) }, \ \ \ \ \ \ \ \ \forall t\geq t_0, \eqno(4.4)$$
where the function $\zeta$ satisfies the assumption (A2).


{\bf Proof.}~~We define the Lyapunov functional
$$L(t):=E(t)+\varepsilon_1 \Psi (t)+\varepsilon_2 \chi (t),\eqno(4.5)$$
where $\varepsilon_i, i=1,2$ are two  positive real numbers which
will be chosen later, and
$$\Psi(t):=\int_\Omega uu_tdx,\eqno (4.6)$$
$$\chi (t):=-\int_\Omega u_t \int_0^t g(t-s)(u(t)-u(s))ds dx. \eqno (4.7)$$

{\emph{\textbf{Remark 4.2}}}~~\emph{We deduce that, for all
$\varepsilon_i>0, i=1,2$, the Lyapunov functional $L(t)$ and the
energy $E(t)$ are equivalent in the sense that there exist two
positive constants $\beta_1,\beta_2$ depending on $\varepsilon_1$
and $\varepsilon_2$ such that
$$\beta_1 E(t)\leq L(t)\leq \beta_2 E(t), \ \ \forall t \geq 0. \eqno(4.8)$$
In fact, using Young's inequality, Poincar\'{e}'s inequality and
(2.3), we have
$$\varepsilon_2 |\chi(t)|\leq
\frac{\varepsilon_2}{2}\|u_t\|_2^2+\frac{\varepsilon_2
(1-l)C_*^2}{2}(g\circ \nabla u)(t)$$
 and
 $$\varepsilon_1 | \Psi (t)|
\leq \frac{\varepsilon_1}{2}\|u_t\|_2^2+\frac{\varepsilon_1
C_*}{2}\|\nabla u\|_2^2.$$ So, we obtain $|L(t)|\leq C E(t)$, where
$C>0$ is a constant. And then (4.8) holds.}

 Now, we will estimate the derivative of $L(t)$
according to the following steps.

\emph{Step 1:}~~Estimate of the derivative of $E(t)$.

Differentiating (4.2), we have
$$E^\prime(t)=\int_\Omega \left( u_tu_{tt}+\left( 1-\int_0^t g(s) ds \right) \nabla u \nabla u_t - \frac{1}{2} g(t) | \nabla u |^2 \right) dx$$
$$+\int_0^t g(t-s)\int_\Omega \nabla u_t (t)\left( \nabla u(t)-\nabla u(s)\right)dxds$$
$$+\frac{1}{2}\int_0^t g^\prime(t-s)\int_\Omega |\nabla u(t)-\nabla u(s)|^2dxds\ \ \ \ \ $$
$$+\frac{\xi}{2}\int_\Omega u_t^2(x,t)dx-\frac{\xi}{2}\int_\Omega e^{-\sigma \tau}u_t^2(x,t-\tau)dx\ \ $$
$$-\frac{\sigma\xi}{2}\int_{t-\tau}^t\int_\Omega e^{-\sigma (t-s)} u_t^2 (x,s) dx ds.\ \ \ \ \ \ \ \ \ \ \ \ \  \eqno(4.9)$$
Using the first equation of (4.1), we obtain
$$E^\prime(t)=\frac{1}{2}\left(g^\prime \circ \nabla u\right)(t)-\frac{1}{2} g(t) \| \nabla u \|_2^2- \mu \int_\Omega u_t(t)u_t(t-\tau)dx \ \ \ \ \ \ \ \ \ \ \ \ \ \ \ \ \ \ \ \ \ \ \ \ \ \ \ $$
$$+\frac{\xi}{2}\|u_t\|_2^2-\frac{\xi}{2}e^{-\sigma \tau}\int_\Omega u_t^2(x,t-\tau)dx-\frac{\sigma\xi}{2}\int_{t-\tau}^t\int_\Omega e^{-\sigma (t-s)} u_t^2 (x,s) dx ds.\eqno(4.10)$$
And then, by Cauchy inequalities, we get
$$E^\prime(t) \leq  \frac{1}{2}\left(g^\prime \circ \nabla u\right)(t)+\left(\frac{|\mu|}{2}+\frac{\xi}{2}\right)\|u_t\|_2^2 +\left(\frac{ |\mu|}{2}-\frac{\xi}{2}e^{-\sigma \tau}\right)\int_\Omega u_t^2(x,t-\tau)dx$$
$$-\frac{1}{2} g(t) \| \nabla u \|_2^2-\frac{\sigma\xi}{2}\int_{t-\tau}^t\int_\Omega e^{-\sigma (t-s)} u_t^2 (x,s) dx ds.\ \ \ \   \eqno(4.11)$$

\emph{Step 2:}~~Estimate of the derivative of $\Psi(t)$.

Differentiating and integrating by parts, we obtain
$$\Psi^\prime(t)=\|u_t\|_2^2+\int_\Omega u\left(\Delta u-\int_0^t g(t-s)\Delta u(s)ds -\mu u_t(t-\tau)\right)dx$$
$$=\|u_t\|_2^2+\int_\Omega \nabla u \int_0^t g(t-s)\left(\nabla u(s)-\nabla u(t)\right)dsdx\ \ \ \ \ $$
$$+\left(\int_0^t g(s)ds-1\right)\|\nabla u\|_2^2-\mu \int_\Omega u(t)u_t(t-\tau)dx\  $$
$$\leq \int_\Omega \nabla u \int_0^t g(t-s)\left(\nabla u(s)-\nabla u(t)\right)dsdx\ \ \ \ \ \ \ \ \ \ \ \ \ \ \ $$
$$+\|u_t\|_2^2-l\|\nabla u\|_2^2-\mu \int_\Omega u(t)u_t(t-\tau)dx.\ \ \ \ \ \ \ \ \ \ \eqno(4.12)$$

By Young's inequality and Lemma 2.1, we get(see [13])
$$\int_\Omega \nabla u \int_0^t g(t-s)\left(\nabla u(s)-\nabla u(t)\right)dsdx\ \ \ \ \ \ \ \ \ \ \ \ \ \ \ $$
$$\leq \delta_1 \|\nabla u\|_2^2+\frac{1}{4\delta_1}\int_\Omega \left(\int_0^t g(t-s)|\nabla u(s)-\nabla u(t)|ds\right)^2 dx$$
$$\leq \delta_1 \|\nabla u\|_2^2+\frac{(1-l)C_*^2}{4\delta_1}\left(g\circ \nabla u\right)(t). \ \ \ \ \ \ \ \ \ \ \  \ (\forall \delta_1>0)\ \ \ \  \eqno(4.13)$$
Also, using Young's and Poincar\'{e}'s inequalities, we have
$$-\mu \int_\Omega u(t)u_t(t-\tau)dx \leq \delta_1 \|\nabla u\|_2^2 + C(\delta_1)\int_\Omega u_t^2(t-\tau)dx. \eqno(4.14)$$
Combining (4.12)-(4.14) and choosing $\delta_1$ small enough, the
estimate
$$\Psi^\prime(t)\leq -\frac{l}{2}\|\nabla u\|_2^2 +C_1\int_\Omega \left(u_t^2(t)+u_t^2(t-\tau)\right)dx+C_2\left(g\circ \nabla u\right)(t)  \eqno(4.15)$$
holds for some positive constants $C_i,i=1,2.$


\emph{Step 3:}~~Estimate of the derivative of $\chi(t)$.

Differentiating (4.7) and integrating by parts, we have
$$\chi^\prime(t)=\left(1-\int_0^t g(s)ds\right)\int_\Omega \nabla u \int_0^t g(t-s)\left(\nabla u(t)-\nabla u(s)\right)dsdx\ \ \ \ \ \ \ \ \ \ \ \ \ \ \ \ \ \ \ \ \ \ \ \ \ \ \ \ \ \ \ \ \ \ \ \ \ \ \ \ \ $$
$$+\int_\Omega \left(  \int_0^t g(t-s)\left(\nabla u(s)-\nabla u(t)\right)ds  \right)^2 dx-\int_\Omega u_t \int_0^t g^\prime (t-s)\left(u(t)-u(s)\right)dsdx$$
$$-\int_0^t g(s)ds \|u_t\|_2^2+\mu\int_\Omega u_t(t-\tau)\int_0^t g(t-s)\left(u(t)-u(s)\right)dsdx.\ \ \ \ \ \ \ \ \ \ \ \ \ \ \ \eqno(4.16)$$
We first estimate the second item of (4.16) as follows:
$$\int_\Omega\left(\int_0^t g(t-s)(\nabla u(s)-\nabla
u(t))ds\right)^2dx\ \ \ \ \ \ \ \ \ \ \ \ \ \ \ \ $$
$$=\int_\Omega\left(\int_0^t \sqrt{g(t-s)}\left(\sqrt{g(t-s)}(\nabla u(s)-\nabla
u(t))\right)ds\right)^2dx$$
$$\leq \int_\Omega\left(\int_0^t g(t-s)\right)\left(g(t-s)(\nabla u(s)-\nabla
u(t))^2 ds\right) dx\ \ \ \ \ \ \ \ $$
$$\leq (1-l)(g\circ \nabla u)(t).\ \ \ \ \ \ \ \ \ \ \ \ \ \ \ \ \ \ \ \ \ \ \ \ \ \ \ \ \ \ \ \ \ \ \ \ \ \ \ \ \ \ \ \ \ \ \ \ \ \ $$
Then, using Young's inequality and Lemma 2.1, we get(see [13,14])
$$\left(1-\int_0^t g(s)ds\right)\int_\Omega \nabla u \int_0^t g(t-s)\left(\nabla u(t)-\nabla u(s)\right)dsdx\ \ \ \ \ \ \ \ \ \ \ \ \ \ $$
$$\leq \delta_2 \|\nabla u\|_2^2+\frac{C_3}{\delta_2}\left(g\circ \nabla u\right)(t),\ \ \ \ \ \ \ \ \ (\forall \delta_2>0)\ \ \ \ \ \  \eqno(4.17)$$
$$-\int_\Omega u_t \int_0^t g^\prime (t-s)\left(u(t)-u(s)\right)dsdx \leq \delta_2\| u_t \|_2^2 - \frac{C_4}{\delta_2}\left(g^\prime\circ \nabla u\right)(t)\eqno(4.18)$$
and
$$\mu\int_\Omega u_t(t-\tau)\int_0^t g(t-s)\left(u(t)-u(s)\right)dsdx \leq \frac{C_5}{\delta_2}\left(g\circ \nabla u\right)(t)+\delta_2\int_\Omega u_t^2(t-\tau)dx, \eqno(4.19)$$
where $C_i,i=3,4,5$ are some positive constants.

Combining (4.16)-(4.19) and (2.3), we obtain
$$\chi^\prime(t) \leq \left(\delta_2-\int_0^t g(s) ds\right) \|u_t\|_2^2 +\delta_2 \|\nabla u\|_2^2+ \frac{C_6}{\delta_2}\left(g\circ \nabla u\right)(t)$$
$$-\frac{C_7}{\delta_2}\left(g^\prime\circ \nabla u\right)(t)+\delta_2\int_\Omega u_t^2(t-\tau)dx, \ \ \ \ \ \ \ \ \ \eqno(4.20)$$
where $C_i,i=6,7$ are some positive constants.

\emph{Step 4:}~~Estimate of the derivative of $L(t)$.

By using (4.3),(4.5),(4.11),(4.15) and (4.20), a series of
computations yields, for $t\geq t_0$,
$$L^\prime(t)\leq \left(\frac{|\mu|}{2}+\frac{\xi}{2}+\varepsilon_1 C_1+\varepsilon_2(\delta_2-g_0)\right)\|u_t\|_2^2+\left(\varepsilon_2\delta_2- \frac{ \varepsilon_1 l}{2}\right)\|\nabla u\|_2^2\ \ \ \ \ \ \ \ \ \ \ \ \ \ \ \ \ \ \ \ \ \ \ \ \ \ \ \ \ \ $$
$$+\left(\frac{1}{2}-\frac{\varepsilon_2 C_7}{\delta_2}\right)\left(g^\prime\circ \nabla u\right)(t)+\left(\varepsilon_1 C_1+\varepsilon_2 \delta_2+\frac{|\mu|}{2}-\frac{\xi}{2e^{\sigma\tau}}\right)\int_\Omega u_t^2(t-\tau)dx$$
$$+\left(\varepsilon_1 C_2+\frac{\varepsilon_2 C_6}{\delta_2}\right)\left(g\circ \nabla u\right)(t)-\frac{\sigma\xi}{2}\int_{t-\tau}^t \int_\Omega e^{-\sigma(t-s)} u_t^2(s)dxds.\ \ \ \ \ \ \ \ \ \ \ \ \ \eqno(4.21)$$

Now, we deduce that, for the positive constants
$\xi,\delta_2,\varepsilon_1,\varepsilon_2$ and $\sigma$, the
following system of inequalities
$$
\left\{
\begin{array}{l}
\frac{|\mu|}{2}+\frac{\xi}{2}+\varepsilon_1 C_1+\varepsilon_2\left(\delta_2-g_0\right)< 0,\nonumber \\
\varepsilon_2 \delta_2-\frac{\varepsilon_1 l}{2}< 0,\nonumber \\
\frac{1}{2}-\frac{\varepsilon_2 C_7}{\delta_2}> 0,\nonumber \\
\varepsilon_1 C_1+\varepsilon_2 \delta_2+\frac{
|\mu|}{2}-\frac{\xi}{2e^{\sigma\tau}}< 0,\nonumber
\end{array}
\right.\eqno(4.22)
$$
is solvable only if we add some suitable conditions to $|\mu|$.

In fact, we can find solutions of (4.22) according to the following
steps.

\emph{Step 1:} We first pick $\delta_2$ small enough such that
$$\delta_2<\min\{\frac{g_0}{4}, \frac{lg_0}{16C_1}\}.\eqno(4.23)$$
Thus, we have
$$\frac{g_0}{8C_1}<\frac{g_0-2\delta_2}{2C_1}.$$

\emph{Step 2:} As long as $\delta_2$ is fixed, we select
$\varepsilon_2$ such that
$$0<\varepsilon_2<\frac{\delta_2}{2C_7}.\eqno(4.24)$$
Then we get
$$\frac{1}{2}-\frac{\varepsilon_2 C_7}{\delta_2}> 0.$$

\emph{Step 3:} Next, we choose $\varepsilon_1$ satisfies the
relation
$$\varepsilon_2 \cdot \frac{g_0}{8C_1}<\varepsilon_1<\varepsilon_2 \cdot \frac{g_0-2\delta_2}{2C_1}.\eqno(4.25)$$
By (4.23) and (4.25), we get
$$\varepsilon_2(g_0-\delta_2)-\varepsilon_1 C_1>\varepsilon_1 C_1+\varepsilon_2 \delta_2>0  \eqno(4.26)$$
and
$$\varepsilon_2 \delta_2-\frac{\varepsilon_1 l}{2}< 0.\eqno(4.27)$$

\emph{Step 4:} Now, we must ensure that the first inequality and the
fourth one in (4.22) hold. That is to say, for the positive
constants $|\mu|,  \sigma, \xi$, the following system of
inequalities
$$
\left\{
\begin{array}{l}
\frac{1}{2}|\mu|+\frac{1}{2}\xi<\varepsilon_2\left(g_0-\delta_2\right)-\varepsilon_1 C_1,\nonumber \\
-\frac{1}{2}|\mu|+\frac{1}{2e^{\sigma\tau}}\xi> \varepsilon_1
 C_1+\varepsilon_2 \delta_2,\nonumber
\end{array}
\right.\eqno(4.28)
$$
must be solvable.

Let $k_1=\varepsilon_2\left(g_0-\delta_2\right)-\varepsilon_1 C_1$
and $k_2=\varepsilon_1
 C_1+\varepsilon_2 \delta_2$. Thus, by (4.26), $k_1$ and $k_2$ are two positive
 constants depending on $g_0$. And then the system (4.28) becomes
$$
\left\{
\begin{array}{l}
|\mu|+\xi<2k_1,\nonumber \\
-|\mu|+\frac{1}{e^{\sigma\tau}}\xi> 2k_2,\nonumber
\end{array}
\right.\eqno(4.29)
$$
By (4.26), we have $k_1>k_2$. Note that $e^{\sigma\tau}\rightarrow
1$ as $\sigma\rightarrow 0$. Thus,  if we choose $\sigma$ small
enough, there exists a positive constant $\xi$ such that
$$2k_2 e^{\sigma\tau}<\xi<2k_1.\eqno(4.30)$$
And then, we have $2k_1-\xi>0$ and
$\frac{\xi}{e^{\sigma\tau}}-2k_2>0$.

 Thus, the system of inequalities (4.29) is
solvable if we choose
$$|\mu|<\min\{2k_1-\xi,\frac{\xi}{e^{\sigma\tau}}-2k_2\}=: a. \eqno(4.31)$$
Here, $a$ is only dependant on $g_0$.

Therefore, (4.22) is solvable under the condition (4.31).

Consequently, from (4.21), there exist two positive constants
$\gamma_1$ and $\gamma_2$ such that
$$\frac{dL(t)}{dt}\leq -\gamma_1E(t)+\gamma_2(g\circ \nabla u)(t),\ \ \ \ \forall t\geq t_0. \eqno(4.32)$$

Similar to the steps in [13], the remaining part of the proof of
inequality (4.4) can be finished. For reader's convenience, we here
write the details as follows:

By multiplying (4.32) by $\zeta$, we arrive at
$$\zeta\frac{dL(t)}{dt}\leq -\gamma_1\zeta E(t)+\gamma_2\zeta  (g\circ \nabla u)(t),\ \ \ \ \forall t\geq t_0. \eqno(4.33)$$
Recalling (A2), (4.2), (4.11) and the first inequality in (4.29), we
get
$$\zeta\frac{dL(t)}{dt}\leq -\gamma_1\zeta E(t)-\gamma_2 (g^\prime \circ \nabla u)(t)\ \ \ \ \ \ \ \ \ \ \ \ \ \ \ \ \ \ \ \ \ \ \ \ \ \ \ \ \ \ \ \ \ \ \ $$
$$\leq -\gamma_1\zeta E(t)-2\gamma_2 E^\prime(t)+2k_1\gamma_2 \| u_t\|_2^2\ \ \ \ \ \ \ \ \ \ \ \ \ \ $$
$$\leq -\gamma_1\zeta E(t)-2\gamma_2 E^\prime(t)+4k_1\gamma_2 E(t),\ \ \ \forall t\geq t_0.\eqno(4.34)$$

Now, we add a restriction condition on $\zeta$, that is, we suppose
that
$$\zeta>\frac{4k_1 \gamma_2}{\gamma_1}.\eqno(4.35)$$
Then, there exists a positive constant $\gamma_3$ such that
$$F^\prime (t)\leq -\gamma_3 \zeta E(t),\eqno(4.36)$$
where $$F(t)=\zeta L(t)+2\gamma_2 E(t)\sim E(t).\eqno(4.37)$$

And then, there exists a positive constant $\gamma_4$ such that
$$F^\prime (t)\leq -\gamma_3 \zeta E(t)\leq -\gamma_4 \zeta F(t).\eqno(4.38)$$
A simple integration of (4.38) over $(t_0,t)$ leads to
$$F(t)\leq F(0)e^{-\gamma_4\zeta(t-t_0)},\ \ \ \forall t\geq t_0.\eqno(4.39)$$
A combination of (4.37) and (4.39) leads to (4.4). The proof of
Theorem 4.1 is thus completed.$\Box$

{\emph{\textbf{Remark 4.3}}}  \emph{We here point out that the
restriction condition (4.35) only makes the admissible space of the
function $g$ a little narrow. But it is not difficult to find that
the constant $\zeta$ can be extended to the case $\zeta=\zeta(t)$,
only if $\zeta(t)$ is a positive nonincreasing differentiable
function with lower bound.}


\bigbreak

 {\bf Open problem}

What results can we hope to get in the case $0<\mu_1<\mu_2$? We
first recall that, in the absence of time-delay, the energy of the
problem (1.2) is dissipative. That is to say, the damping item is a
``good" item for the energy decay. So it is not difficult to imagine
that, if $|\mu_1|<\mu_2<a$(Here, $a$ is defined by (4.31)), we also
can obtain the energy decay result. But, if $a<|\mu_1|<\mu_2$, it
could be instability. And this problem is still open.


\bigbreak

 {\bf Acknowledgments}

The authors would like to express their sincere gratitude to the
anonymous referee for his/her valuable comments and useful
suggestions on the manuscript of this work.

\bigbreak
\def\item{\par\hang\textindent}
\vskip 10pt{\small \centerline{\Large\bf References} \medbreak

\item{[1]}
        J. Hale and S. Verduyn Lunel, ``Introduction to Functional Differential Equations", Volume 99 of ``Applied Mathematical Sciences",
        {\sl Springer-Verlag, New York.}, 1993.
 \item{[2]}
        R. Datko, J. Lagnese, M.P. Polis, An example on the effect of time delays in boundary feedback stabilization of wave equations,
        {\sl SIAM J.Control Optim}, 1986, {\bf 24}:152-156.
 \item{[3]}
        R. Datko, Not all feedback stabilized hyperbolic systems are robust with respect to small time delays in their feedbacks,
        {\sl SIAM J.Control Optim}, 1988, {\bf 26}:687-713.
 \item{[4]}
         R. Datko, Two examples of ill-posedness with respect to time delays revisited methods,
         {\sl IEEE Trans. Automat. Control}, 1997,{\bf 42}:511-515.
\item{[5]}
        G.Q. Xu, S.P. Yung and L.K. Li, Stabilization of wave systems with input delay in the boundary control,
         {\sl ESAIM Control Optim. Calc.Var.}, 2006,{\bf 12(4)}:770-785.
\item{[6]}
        S. Nicaise, C. Pignotti, Stability and instability results of the wave equations with a delay term in the boundary or internal feedbacks,
         {\sl SIAM J.Control Optim}, 2006, {\bf 45}:1561-1585(electronic).
\item{[7]}
        S. Nicaise, J. Valein, Stabilization of the  wave equation on
        1-D networks with a delay term in the nodal feedbacks,
         {\sl Netw. Heterog. Media}, 2007,{\bf 2(3)}:425-479(electronic).
\item{[8]}
        S. Nicaise, C. Pignotti, Stabilization of the wave equation
        with boundary or internal distributed delay, {\sl Diff. Int. Equs.}, 2008,{\bf 21(9-10)}:935-958.
\item{[9]}
        M.M. Cavalcanti, V.N. Domingos Cavalcanti, J.A. Soriano,
        Exponential decay for the solution of semilinear
        viscoelastic wave equations with localized damping, {\sl E. J. Differ. Int. Equ.}, 2002,{\bf 44}:1-14.
\item{[10]}
        S. Berrimi, S.A. Messaoudi, Existence and decay of solutions of a viscoelastic equation with a nonlinear source,
        {\sl Nonl. Anal.}, 2006, {\bf 64}:2314-2331.
\item{[11]}
        F. Alaban-Boussouira, P. Cannarsa, D. Sforza, Decay
        estimates for second order evolution equations with memory.
        {\sl J. Funct. Anal.}, 2008,{\bf 254(5)}:1342-1372.
\item{[12]}
         M.M. Cavalcanti, H.P. Oquendo, Frictional versus viscoelastic damping in a semilinear wave equation,
        {\sl SIAM J.Control Optim.}, 2003,{\bf 42(4)}:1310-1324.
\item{[13]}
         M. Kirane, B. Said-Houari, Existence and asymptotic stability of a
        viscoelastic wave equation with a delay,
        {\sl Z. Angew. Math.Phys.}, 2011,{\bf 62}:1065-1082.
\item{[14]}
        W.J. Liu, General decay of the solution for viscoelastic wave
        equation with a time-varying delay term in the internal
        feedback, {\sl arXiv: 1208.4007v1}, 20 Aug, 2012.
\item{[15]}
        J.E. Munoz Rivera, M. Naso, E. Vuk, Asymptotic behavior
        of the energy for electromagnetic system with memory,
        {\sl J. Math. Meth. Appl. Sci.}, 2004, {\bf 25(7)}:819-941.
\item{[16]}
        S.A. Messaoudi, General decay of solutions of a viscoelastic
        equation, {\sl J. Math. Anal. Appl.}, 2008, {\bf 341}:1457-1467.
\item{[17]}
         J.-L. Lions, Quelques m\'{e}thodes de r\'{e}solution des
        probl\`{e}mes aux limites non lin\'{e}aires. {\sl Dunod},
        1969.
\item{[18]}
        N.-e. Tatar, Arbitrary decay in linear viscoelasticity, {\sl J. Math. Phys.}, 2011, {\bf 52}: 01350201-01350212.
\item{[19]}
        M.I. Mustafa, S.A. Messaoudi, General stability result for
        viscoelastic wave equations, {\sl J. Math. Phys.}, 2012, {\bf 53}: 05370201-05370214.
\item{[20]}
        M. Gugat, Boundary feedback stabilization by time delay for one-dimensional wave
        equations, {\sl IMA J Math Control Info}, 2010, {\bf 27(2)}:189-203.
\item{[21]}
        M. Gugat, G. Leugering, Feedback stabilization of quasilinear hyperbolic systems with varying
        delays, {\sl Methods and Models in Automation and Robotics(MMAR), 2012 17th International Conference on}, 27-30 Aug. 2012:125-130.
\item{[22]}
        C. Pignotti, A note on stabilization of locally damped wave equations with time
        delay, {\sl Systems \& Control Letters}, 2012, {\bf 61(1)}:92-97.
\item{[23]}
        L. C. Evans, Partial Differential Equations, Graduate Studies in Mathematics, Volume 19,{\sl American Mathematical Society}, 1997.

\end{document}